\documentclass[10pt]{article}
\usepackage{euscript}
\title{\bf \LARGE Vanishing of the 
$\tilde K{\cN}il$ groups : localization methods}
\author{\bf Frank BIHLER}
\date{\empty}
\usepackage{fancyheadings,ifthen,pifont}
\input xy
\xyoption{all}
\usepackage{amssymb,amsfonts,amsmath,mathrsfs,stmaryrd,bbm}
\usepackage{color,epsfig,graphicx,floatflt,footmisc}

\newtheorem{defin}{Definition \\}
\newtheorem{theo}{Theorem \\}
\newtheorem{prop}{Proposition \\}
\newtheorem{lemme}{Lemma \\}

\def\no{\noindent}

\newfam\msyfam

\def \Z {{\mathbbm{Z}}}

\def \Q {{\mathbbm{Q}}}
\def \A {{\frak{A}}}
\def \B {{\frak{B}}}
\def \C {{\frak{C}}}
\def \L {{\frak{L}}}
\def \N {{\mathbbm{N}}}
\def \T {{\mathbbm{T}}}
\def \F {{\mathbbm{F}}}
\def \W {{\frak{W}}}

\font\tgoth=eufm10
\font\sgoth=eufm7
\font\ssgoth=eufm5
\newfam\gothfam
\textfont\gothfam=\tgoth
\scriptfont\gothfam=\sgoth
\scriptscriptfont\gothfam=\ssgoth

\def \bI {{\mathbbm{1}}}

\def \cH {{\EuScript{H}}}
\def \cD {{\mathscr{D}}}
\def \cA {{\mathscr{A}}}
\def \cB {{\mathscr{B}}}

\def \cN {{\EuScript{N}}}

\font\tcal=eusm10
\font\scal=eusm8
\newfam\calfam
\textfont\calfam=\tcal
\scriptfont\calfam=\scal

\font\callfont=rsfs12
\newfam\callfam
\textfont\callfam=\callfont
\def \call{\fam\callfam\callfont}

\def\build#1_#2^#3{\mathrel{\mathop{\kern 0pt#1}
\limits_{#2}^{#3}}}
\newdir{ >}{{}*!/-8pt/@{>}}

\begin{document}

\maketitle

\vspace{1cm}
\no
\begin{abs}

\no
The aim of this article is to introduce Vogel's 
localization theorem for classes of 
${\cD}$-complexes : this generalization of 
Waldhausen's localization theorem is especially 
useful and powerful in that it gives an explicit
and computable description of the local objects.
Next we present an excision theorem for transverse
classes. Finally, we apply these methods to deduce
partial results on Vogel's [Conjecture] on a 
regular ring $R$ ( cf \cite{bil1} ). 
\end{abs}

\vspace{1cm}
\no
\begin{key}

\no
COMPLEX OF DIAGRAMS -- LOCALIZATION -- EXCISION \\
VANISHING OF THE $\tilde K{\cN}il$ GROUPS -- 
VOGEL'S CONJECTURE \\
\end{key}

\tableofcontents
\thispagestyle{empty}

\newpage
\section{Main Theorem}
\no
We shall expose at first the localization 
theorem, powerful tool introduced by Vogel 
in \cite{vog2}, as a motivation to develop
further notations and constructions, and then 
apply it to our [Conjecture] on a regular ring 
$R$ : namely, give a description of the algebraic
K-theory spectrum $\tilde K{\cN}il(R;S)$ for 
$S$ a $R$-bimodule, flat from the left, reduce 
it through explicit localization computations, 
and ideally prove the vanishing of these 
obstruction groups to excision. In a few words, 
the reader can think about $\cD$-complexes as 
representations of a diagram $\cD$ in the category 
of chain complexes : h-small and h-compact are 
finiteness hypothesis to insure the existence of 
algebraic K-theory spectra, and the 'exactness' 
condition reminds about Waldhausen's 'stable under
extensions subcategory' or about Verdier's 'thick 
subcategory' needed for a good localization. These 
form the natural setting for a generalization of the 
different categories of $\cN il$ introduced by Waldhausen
in \cite{wald2}. Precise definitions will just follow. 
     
\begin{theo}{\em [Localization] \cite{vog2}}\\
Let $\A  \subset \B $ be two exact 
classes of ${\cD}$-complexes. Suppose that 
the class $\B $ is h-small and h-compact, 
and that the class $\A $
is stable in $\B $ under direct summands.
Let $\L$ be the class of all 
$\A $-local ${\cD}$-complexes $L$, such 
that there exists $X \in \B $ and an 
$\tilde\A $-equivalence $X \to L$. 
Then $\L$ is an exact class, and we have 
an homotopy equivalence of spectra :
$K(\B,\A) \simeq K(\L)$.
\end{theo}
\vspace*{-.5cm}
\section{${\cD}$-complexes}
\no
We shall now expose, by order of increasing 
abstraction, the numerous definitions required 
to understand the hypothesis of the above 
localization theorem.
\begin{defin}.\\
\ding{172} A ``{\bf diagram of bimodules}'' 
$\cD$ is given by :
\\$\bullet$ an oriented graph $\Gamma$ \\
$\bullet$ for each vertex $a \in \Gamma_0$, an
associated ring $A_a$ \\ $\bullet$ for each edge
$(f:a \to b) \in \Gamma_1$, an associated bimodule
$S_f$ with \\ \hspace*{.2cm} an action of $A_b$ on 
the left, and an action of $A_a$ on the right.\\
Consider now the diagram $\cD$ fixed. \\
\ding{173} A ``{\bf $\cD$-module}'' 
is given by :\\
$\bullet$ for each vertex $a$, an associated 
right $A_a$-module $M_a$ \\
$\bullet$ for each edge $(f:a \to b)$, an 
associated map of modules $M_a \to M_b \otimes 
S_f$ \\
The $\cD$-modules form an abelian category.\\ 
We can thus take its derived category and 
define :\\
\ding{174} A ``{\bf $\cD$-complex}'' is a complex 
of $\cD$-modules. \\
\ding{175} A $\cD$-module is called 
``{\bf nilpotent}'' if there exists a length 
$n$ such that every \\ \hspace*{.2cm} 
composition of more than $n$ linear maps vanishes. 
Similarly, a $\cD$-complex is \\
\hspace*{.2cm} called ``{\bf h-nilpotent}''
or homotopy-nilpotent, or nilpotent up to homotopy,
if \\ \hspace*{.2cm} there exists a length $n$ such 
that every composition of more than $n$ chain maps\\
\hspace*{.2cm} is null-homotopic.  
\end{defin}
\no
{\bf Convention :} \\
As the aim of this article is to treat objects fit 
to algebraic K-theory, we shall from now on take 
this convention : a {\bf ``ring''} $A$ shall mean a unitary 
associative ring; a {\bf ``module''} $M$ shall mean a right 
projective module (unless otherwise stated); a {\bf ``bimodule''}
$S$ shall mean a flat on the left bimodule; a {\bf ``complex''} 
shall mean a $\Z$-graded module $M_*$ equipped with a 
differential $\partial:M_n \to M_{n-1}$ such that 
$\partial \circ \partial =0$. With this convention, we only ask 
a $\cD$-complex $C_*$ to be composed by complexes $C_a$ 
on each vertex $a$ of the graph $\Gamma$, projective in each degree.
Now a {\bf ``chain map''} shall mean a morphism of complexes. 
Define a {\bf ``finite complex''} $C_*$ to be a complex such that 
$\oplus C_n$ is finitely generated. Finally, a {\bf ``h-finite 
complex''} $D_*$ shall mean a complex equipped with a finite
complex $C_*$ and two chain maps $f:C_* \to D_*$ 
and $g:D_* \to C_*$, such that $f \circ g \sim 
Id:D_* \to D_*$ and $g \circ f \sim Id:C_* \to C_*$. 
We shall denote $\C_R$ the class of h-finite complexes 
on the base ring $R$. \\
  
\no
{\bf Examples :} \\
Friedhelm Waldhausen introduced three categories 
$\cN il(R;S)$, $\cN il(R;S,T)$ and $\cN il(R;S,T,U,V)$ 
to study the algebraic K-theory of a tensor algebra, a free 
product of rings, and respectively a Laurent extension 
(or HNN-extension) of rings; these categories can be naturally 
viewed as categories of nilpotent $\cD$-modules, on the diagrams 
$\cD_1,\cD_2$ and $\cD_3$ hereafter. In each case, the involved 
category of $\cD$-modules gives the obstructions to excision 
detailed in \cite{wald2}. \\
$\bullet$ The category ${\cN}il(R;S)$ defined for 
the {\bf tensor algebra $R[S]$} can be represented 
by the following diagram $\cD_1$ of bimodules;
similarly, one can draw the 'picture' of a 
$\cD_1$-module, a nilpotent $\cD_1$-module, 
a $\cD_1$-complex, and a homotopy-nilpotent 
$\cD_1$-complex :  \\
$$\xymatrix@-1.5pc{
\cD_1& =& \bullet_R \ar@(dr,ur)[]_S \\
\cD_1 \text{-module}&=& \bullet_{M_R} 
\ar@(dr,ur)[]_f &
(\text{with}\;f:M_R \to M_R \otimes_R S 
\;\;\text{linear} ) \\
\text{nilpotent} \cD_1 \text{-module}&=& 
\bullet_{M_R} \ar@(dr,ur)[]_f &
( \exists n \ge 0, 
f^n=0 : M_R \to M_R \otimes_R S^n ) \\
\cD_1\text{-complex}&=& \bullet_{C_R} 
\ar@(dr,ur)[]_f&
\hskip.5cm ( \text{with}\;f:C_R \to C_R \otimes_R S 
\hskip.2cm \text{chain map} )\\
\text{h-nilpotent} \cD_1\text{-complex}&=&
\bullet_{C_R} \ar@(dr,ur)[]_f & 
( \exists n \ge 0, 
f^n \sim 0 : C_R \to C_R \otimes_R S^n )}$$
$\bullet$ Let's consider ${\cN}il(R;S,T)$ defined
for the {\bf generalized free product} of rings,
where the two pure embeddings $\alpha,\beta$ admit 
splittings as $R$-bimodules. \vspace*{-.3cm}
$$\xymatrix@-1pc{
R \ar@{ >->}[r]^{\alpha} \ar@{ >->}[d]_{\beta} & 
R \oplus S \ar@{ >.>}[d] \\
R \oplus T \ar@{ >.>}[r] & \text{pushout} 
}$$ 
Here the naturally associated diagram $\cD_2$
of bimodules is of the shape : \vspace*{-.3cm}
$$\xymatrix@-.5pc{ 
\cD_2 &=& \bullet_R \ar@/^/[r]^S & 
\bullet_R \ar@/^/[l]^T \\
\text{h-nilpotent} \cD_2 \text{-complex} &=&
\bullet_{C_R} \ar@/^/[r]^f &
\bullet_{D_R} \ar@/^/[l]^g}$$ \vspace*{-.3cm}
$$\parbox{10.6cm}{$( \text{with} \; f:C_R \to D_R 
\otimes_R S \; \text{and} \; g:D_R \to C_R 
\otimes_R T \;\;\text{chain maps}\\ \text{and} 
\;\exists n \ge 0, \;\text{such that}\; 
(g \circ f)^n \sim 0 :
C_R \to C_R \otimes_R (S \otimes_R T)^n )$}$$
$\bullet$ Let's consider ${\cN}il(R;S,T,U,V)$ 
defined for the {\bf Laurent extension} of rings, where 
the two pure embeddings $\alpha,\beta$ admit 
splittings as $R$-bimodules.
$$\xymatrix{
(R \ar@<.5ex>@{ >->}[r]^\alpha 
\ar@<-.5ex>@{ >->}[r]_\beta &
A \ar@{.>}[r] & \text{HNN-ext} &
\text{with} \;A=S \oplus U =T \oplus V; S=\alpha R;
T=\beta R)}$$
Here the naturally associated diagram $\cD_3$ is of
the shape : \vspace*{-.3cm}
$$\xymatrix@-.5pc{ 
\cD_3 &=& \bullet_R \ar@/^/[r]^U \ar@(dl,ul)^S & 
\bullet_R \ar@/^/[l]^V \ar@(dr,ur)_T\\
\text{h-nilpotent} \cD_3 \text{-complex} &=&
\bullet_{C_R} \ar@/^/[r]^f \ar@(dl,ul)^h &
\bullet_{D_R} \ar@/^/[l]^g \ar@(dr,ur)_i}$$ 
\vspace*{-.3cm}
$$\parbox{8.5cm}{$( \text{with} \; f:C_R \to D_R 
\otimes_R U \; \text{and} \; g:D_R \to C_R 
\otimes_R V \; \text{and} \\ h:C_R \to C_R 
\otimes_R S \;\text{and} \; i:D_R \to D_R
\otimes_R T \;\text{chain maps}\\ \text{and} 
\;\exists n \ge 0, \; \text{such that 
every composed map of length} \\
n \;\text{is null-homotopic; in particular this  
composed map :} 
\\ h^{n-4}\circ (g \circ f)^2 \sim 0 : 
C_R \to C_R \otimes_R (U \otimes_R V)^2
\otimes_R S^{n-4} )$}$$ 
In terms of algebraic K-theory, we have the 
following decomposition :
$$K \cN il(R;S) = K(R) \oplus \tilde K \cN il(R;S)$$
More generally, for each diagram of bimodules $\cD$,
we can use the notion of h-nilpotent $\cD$-complexes
( here we use h-finite complexes ) to define a 
K-theory spectrum $\tilde K \cN il(\cD)$ that will extend 
the usual notions presented above, and first 
introduced by Waldhausen in \cite{wald2}.
\section{Exact Classes of $\cD$-complexes}
\begin{defin}.\\
Let $\cD$ be a fixed diagram of bimodules, 
and $\A$ a class of $\cD$-complexes.\\
\ding{176} The class $\A$ is called 
``{\bf exact}'' if :\\
$\bullet$ $\A$ contains all finite acyclic 
$\cD$-complexes (ie acyclic on each vertex of $\Gamma$)\\
$\bullet$ $\A$ verifies the '2/3 axiom' :
let $0 \to A \to B \to C \to 0$ be a short 
exact \\\hspace*{.2cm}sequence of $\cD$-complexes, 
if two of them are in $\A$, then so is the 
third. \\ 
\ding{177} The ``{\bf completion}'' of $\A $, 
denoted by $\tilde \A $, is the smallest 
exact class \\ \hspace*{.2cm}containing 
$\A $ and stable under direct sum : 
$\forall i \in I, A_i \in \tilde\A  
\Rightarrow \oplus A_i \in \tilde\A $.\\
\ding{178} A morphism $f$ of $\cD$-complexes 
is a ``{\bf $\A $-equivalence}'' if 
the usual \\ \hspace*{.2cm}mapping cone 
$C(f) \in \A $. 
\end{defin}
Recall the classical construction issued 
from Waldhausen's cylinder-functor $\T$ :
\vspace{-0.95cm}
$$\xymatrix{A \ar@{ >->}[r]^{j_1} \ar[dr]_f &
T(f) \ar@{->>}[d]^p & B \ar@{ >->}[l]_{j_2} 
\ar@{=}[dl] \\ &B} \hspace{2cm}  
\parbox{5.2cm}{\vspace*{1cm} 
\no $\xymatrix@-.6pc{0 \to A \ar@{ >->}^{j_1}[r] & 
T(f) \ar@{->>}[r] &C(f) \to 0}$  \\ So 
$C(f) \in \A $ by the '2/3 axiom'.}$$
\section{Finiteness Properties}
\begin{defin}.\\
Let $\cD$ be a fixed diagram of bimodules, and 
$\A $ a class of $\cD$-complexes. \\
\ding{179} The class $\A $ is 
``{\bf h-small}'' \footnote[2]{This condition avoids set-theoretic 
problems for the existence of the K-theory $K(\A )$.} 
or homotopy-small, or small up to homotopy if \vspace{-.3cm}
\begin{center}
$\exists (X_i)_{i \in I},\forall Y \in 
\A , \exists i, \exists \xymatrix@-1pc{X_i 
\ar[r]^{\simeq} & Y}$
\end{center} \vspace{-.3cm}
\hspace*{.2cm}Here $(X_i)_{i \in I}$denotes a 
\underline{set} of $\cD$-complexes, and $\simeq$ a
homology equivalence.\\
\ding{180} The class $\A $ is 
``{\bf h-compact}'' or homotopy-compact, or compact \\
up to homotopy if 
$\forall X \in \A , \forall (Y_i)_{i \in I} \in \A ^I,
\forall (f:X \to \oplus Y_i) \in Fl(\tilde\A )$, \\ 
\vspace*{-.4cm} \\
\no
$\exists \xymatrix@-1pc{
(g:X' \ar[r]^(.55){\simeq} & X)\in} Fl(\A )$
such that $f \circ g$ factorizes through a finite 
sum.\\
\ding{181} A $\cD$-complex $X$ is 
``{\bf $\A $-local}'' if 
every $\cD$-complex morphism 
$f:Y \to X$ \\ \hspace*{.2cm}with $Y \in \A $ 
factorizes through an acyclic $\cD$-complex.  
\end{defin}
\section{Easy Applications}
{\bf $1^{st}$ example :}\\
Let $\B $ be the class of usual h-finite 
$\Z$-complexes \mbox{( here the diagram $\cD_4$ is 
reduced} to one vertex with the ring $\Z$ ). 
And consider the sub-class $\A  \subset 
\B $ of $\cD_4$-complexes that are 
$\Q$-acyclic. Then $\A,\B$ are exact classes 
of $\cD_4$-complexes. The class $\B $ is 
h-small and h-compact, and $\A $ is stable in 
$\B $ under direct summands. Now Theorem 1 
gives us the class $\L$ of all $\Z$-complexes
with the rational homotopy type of a finite complex,
and an homotopy equivalence of spectra : 
$K(\B,\A) \simeq K(\L)$. Here we 
consider $\L'$ the class of all finite 
$\Q$-complex. The natural inclusion 
$\L' \subset \L$ 
induces an equivalence 
on the algebraic K-theory spectra by Waldhausen's
approximation theorem, and $K(\L)=
K(\L')=K(\Q)$. We can also identify :
$K(\B )=K(\Z)$ and $K(\A )=\oplus 
K(\F_p)$ ( the direct sum here is over the prime 
numbers $p$ ) by Quillen's D\'evissage. Hence we 
get the classical commutative \mbox{localization 
long exact sequence, usually obtained via the Ore 
condition.}   
$$ \ldots \to \oplus K_i(\F_p) \to K_i(\Z) 
\to K_i(\Q) \to \oplus K_{i-1}(\F_p) \to \ldots $$ 

\no
{\bf $2^{nd}$ example :}\\
Consider now the following diagram of bimodules :
\vspace*{-.3cm}
$$\xymatrix{ 
\cD_5 &=& \bullet_A \ar@/^/[r]^S & 
\bullet_B \ar@/^/[l]^T}$$ \\ \vspace*{-.8cm} \\
We note $\B $ the class of h-nilpotent 
$\cD_5$-complexes, ie $\B =\{ (C_A,C_B;\alpha,
\beta) \}$ with $C_A,C_B$ h-finite complexes, 
$\alpha:C_A \to C_B \otimes S$ and 
$\beta:C_B \to C_A \otimes T$ verifying a nilpotency
condition : $\exists n \ge 0, (\alpha \beta)^n 
\sim 0$. Consider now the sub-class $\A  
\subset \B $ of $\cD_5$-complexes such that 
$C_A$ is contractible. Then the class $\B $ 
is h-small and h-compact, and the class $\A $
is stable in $\B $ under direct summands. 
Now look at the $\A $-local objects : 
we test on $(0,C_B';0,0) \in \A $ to obtain 
the following data :
$$\xymatrix@-1pc{
0\ar@<.5ex>[r] \ar[d] &C_B' \ar@<.5ex>[l] \ar[d]^f\\
C_A\ar@<.5ex>[r]^\alpha &C_B \ar@<.5ex>[l]^\beta}$$ 
Every vertical map of $\cD_5$-complex is thus given 
by $f:C_B' \to Ker \beta$. We can suppose that 
$\beta$ is surjective, so $\A $-local implies
that $\beta$ is a homology equivalence. Conversely,
every $\cD_5$-complex with $\beta$ a homology 
equivalence is easily seen to be $\A $-local. 
Here the condition of $\tilde \A $-equivalence
says that in fact $C_B$ is h-finite. 
Now, such a data is just $(C_B,\alpha \beta)$, that
is \footnotemark[4] an element of 
$\cN il(B;T \otimes_A S)$, thus
Waldhausen's approximation theorem gives the 
equivalence of spectra : $K(\L) \simeq 
K \cN il(B;T \otimes_A S)$. Now Theorem 1 gives 
this homotopy fibration : $\xymatrix@-1pc{
K(A) \ar@{ >->}[r] & K(\B ) 
\ar@{->>}[r] & K \cN il(B;T \otimes_A S)}$ 
symmetrically, we have also this other \vskip-.15cm
\no homotopy fibration :
$\xymatrix@-1pc{
K(B) \ar@{ >->}[r] & K(\B ) 
\ar@{->>}[r] & K \cN il(A;S \otimes_B T)}$. Hence 
comparing the two lines above gives the 
identification : $K(\B ) \simeq K(A) \oplus 
K(B) \oplus$ ``defect''. 
\begin{prop}.\\
Let $A,B$ be two rings, and let ${}_A S_B$ and ${}_B T_A$ be two 
bimodules, flat on the left. Then we have a homotopy equivalence : 
$\tilde K \cN il(B;T \otimes_A S) \simeq 
\tilde K \cN il(A;S \otimes_B T)$.
\end{prop}
\footnotetext[4]{The equivalence of the different 
definitions of $\cN il(R;S)$ is proven in \cite{bil2}.}
\no
{\bf Use :}\\
If $S$ is flat from the left and $A$ Noetherian 
regular, then \cite{wald2} implies that 
$\tilde K \cN il(A;S \otimes_B T)$ is contractible.
Thus so is the spectrum $\tilde K \cN il(B;T 
\otimes_A S)$ ! Hence some tricky reductions of 
the $\cN il$ groups on diagrams ... \\

\no
{\bf $3^{rd}$ example :} \\
If we keep the notations of \cite{bil2}, for a 
regular ring $R$ and a flat from the left bimodule 
$S$, we have the following homotopy equivalences
of non-connective spectra :
$\tilde K \cN il(R;S) \simeq K(\cB,\cA) = 
\varinjlim K(\cB_n,\cA_n)$. So an approach of 
Vogel's [Conjecture] based on localization should 
study first the graded case for $\cB$ :  
($\B =\cB_{n+1}$, $\A =\cB_n$), and 
then the graded case for \mbox{$\cA$ :
($\B =\cA_{n+1}$, $\A =\cA_n$).} 
It appears that the $\cB$ case breaks out easily, 
but the nilpotent case $\cA$ gives only partial 
information, for lack of a good interpretation of 
local objects. 
\mbox{We postpone the complete description of these 
computations till section 8.}

\section{Proof of the Main Theorem}
\no 
This rather technical proof shall be 
decomposed in eight elementary lemmas. 
Notations are those of Theorem 1.

\begin{lemme}.\\
A class $\A $ of $\cD$-complexes is exact, 
if and only if $\A $ is stable under the 
$Coker$ of a cofibration, $\A $ is stable 
under desuspension, and $\A $ contains the 
acyclics.
\end{lemme}
\begin{demo}
Let $0 \to X \to Y \to Z \to 0$ be a short exact 
sequence of $\cD$-complexes. Suppose that $Z$ is 
in $\A $. There exists an epimorphism 
$E \to Y$ with $E$ acyclic. We can now ``make the 
triangle turn'' : if we note $s^{-1}Z$ the 
desuspension of $Z$ in $\A $, then 
$0 \to s^{-1}Z \to X \oplus E \to Y \to 0$ is also 
a short \mbox{exact sequence.} Now, if $X$ is in 
$\A $, this shows that $Y$ is also in 
$\A $. Hence the class $\A $ is stable 
under extensions. Conversely, if $Y$ is in 
$\A $, then $X \oplus E$ is also, hence 
$X=Coker[E \to X \oplus E]$ is \mbox{in $\A $ 
too. Thus the class $\A $ verifies the 
'2/3 axiom', $\A $ contains the acyclics, so 
$\A $ is exact. $\blacksquare$} 
\end{demo}

\begin{lemme}.\\
Let $\A $ be an exact class of 
$\cD$-complexes. \\Let $\L\A$ be the class
of $\A $-local $\cD$-complexes. Then 
the class $\L\A$ is also exact. 
\end{lemme}
\begin{demo}
The class $\L\A$ contains the acyclics. Like 
$\A $, the class $\L\A$ is stable under 
desuspension. So we need to prove that $\L\A$ 
is stable under the $Coker$ of a cofibration. Let 
$0 \to X \to Y \to Z \to 0$ be a short exact 
sequence, with $X$ and $Y$ in $\L\A$. Consider
a $\cD$-complex $U$ in $\A $ and a 
$\cD$-complex morphism $f:U \to Z$. 
There exists an epimorphism $E \to U \prod_Z Y$ with
$E$ acyclic and onto the pullback $U \prod_Z Y$. 
Note $V=Ker[E \to U]$ : the exact class $\A $
contains $U$ and $E$, thus $V$ too. We obtain the 
following commutative diagram :
$$\xymatrix@-1pc{
0 \ar[r] & X \ar[r] & Y \ar[r] & Z \ar[r] & 0 \\
0 \ar[r] & V \ar[u] \ar[r] & E \ar[r] \ar[u] &
U \ar[r] \ar[u]_f & 0}$$
The map $V \to X$ admits a factorization through an 
acyclic $\cD$-complex $F$. Build the pushout 
$W=F \coprod_V E$ : it's an object in $\A $, 
thus the map $W \to Y$ admits a factorization 
through an acyclic $\cD$-complex $G$. We can suppose
that $F \to G$ is a cofibration ( if need be, we can
add an acyclic to $G$ ), and denote $H=Coker[F \to 
G]$. Then the map $f$ admits a factorization through
$H$ acyclic ( cf the following diagram ), and 
the class $\L\A$ is exact. $\blacksquare$ 
$$\xymatrix@-1pc{
0 \ar[r] & X \ar[r] & Y \ar[r] & Z \ar[r] & 0 \\
0 \ar[r] & F \ar[r] \ar[u] & G \ar[r] \ar[u] &
H \ar[r] \ar[u] & 0 \\
0 \ar[r] & V \ar[u] \ar[r] & E \ar[r] \ar[u] &
U \ar[r] \ar[u] & 0}$$ 
\end{demo} 

\begin{lemme}.\label{trois}\\
Let $\A $ be an exact, h-small and h-compact 
class of $\cD$-complexes. Let $X$ be any 
$\cD$-complex. Then there exists an 
$\A $-local $\cD$-complex $Y$ and a 
$\tilde \A $-equivalence $X \to Y$.
\end{lemme}
\begin{demo}
We shall build by induction a series $X=Z_0 \to Z_1
\to Z_2 \to \ldots$ of $\cD$-complexes, where all 
maps $Z_i \to Z_{i+1}$ are cofibrations and $\tilde
\A $-equivalences, and where the colimit $Z$ 
shall be a $\A $-local $\cD$-complex. As the 
class $\A $ is h-small, denote $\{X_\lambda\}
_{\lambda \in \Lambda}$ the set of $\cD$-complexes 
that generate $\A $. Suppose that $Z_i$ is 
built for every $0 \leq i \leq n$. Let $T_n$ be the 
set of pairs $(\lambda,u)$ with $\lambda \in 
\Lambda$ and $u:X_\lambda \to Z_n$. Denote $U_n$ the
direct sum of all $X_\lambda$ for $(\lambda,u) 
\in T_n$. We have a canonical map from $U_n$ to 
$Z_n$. Let $U_n \to E_n$ be a cofibration into an 
acyclic $E_n$, and build the pushout $Z_{n+1}=Z_n 
\coprod_{U_n} E_n$. Then the map $Z_n \to Z_{n+1}$ 
is a cofibration; it's a $\tilde 
\A $-equivalence since $X_\lambda \in 
\A $ implies $U_n \in \tilde\A $. Let 
now $Z$ be the colimit of this series, and $f:Y 
\to Z$ any map from any $Y$ in $\A $. Let 
$\alpha :\oplus Z_n \to \oplus Z_n$ be the 
difference between the identity and the 
stabilization map $Z_n \to Z_{n+1}$. There exists 
an epimorphism $E \to \oplus Z_n \prod_Z Y$ from 
an acyclic $E$ onto the pullback $\oplus Z_n 
\prod_Z Y$. We obtain the following commutative 
diagram : \vspace*{-.2cm}
$$\xymatrix@-.6pc{
0 \ar[r] & \oplus Z_n \ar[r]^\alpha & 
\oplus Z_n \ar[r] & Z \ar[r] & 0 \\
0 \ar[r] & K \ar[u] \ar[r] & E \ar[u] \ar[r] &
Y \ar[u]_f \ar[r] & 0}$$
Now we shall use twice the hypothesis of h-compacity
on $\A $ to reduce ( modulo homology 
equivalence ) the two infinite sums above to finite 
sums : there exists a homology equivalence 
$K' \to K$ such that the composed map 
$K' \to K \to \build{\oplus}_0^\infty Z_n$ 
admits a factorization through $\build{\oplus}
_0^p Z_n$. We can suppose that the map $K' \to K$ 
is a cofibration ( if need be, we can add an acyclic
to $K$ and $E$ ). Now the map $\alpha$ sends $\build
{\oplus}_0^p Z_n$ to $\build{\oplus}_{n+1}^\infty 
Z_n$; by quotient, the middle vertical map sends 
$E/K'$ to $\build{\oplus}_{p+2}^\infty$. By 
h-compacity, there exists a homology equivalence 
$V \to E/K'$ that sends $V$ on a finite sum. 
Build now the pullback $E'=E \coprod_{E/K'} V$. 
We have thus built a finite factorization of the 
diagram above, for a certain finite integer $q$ 
( the cokernel $Y'$ is homologically equivalent to 
$Y$ ) : \vspace*{-.2cm} 
$$\xymatrix@-.6pc{
0 \ar[r] & \build{\oplus}_0^{q-1} Z_n \ar[r]^\alpha 
& \build{\oplus}_0^q Z_n \ar[r] & Z_q \ar[r] & 0 &
Z_q \ar[r] & Z \\
0 \ar[r] & K' \ar[r] \ar[u] & E' \ar[r] \ar[u] & 
Y' \ar[r] \ar[u]_f & 0 & Y' \ar[r] \ar[u] & 
Y \ar[u]_f }$$
Now $Y'$ being in $\A $ is equivalent to a 
certain $X_\lambda$, then the map $u:X_\lambda \to 
Y' \to Z_q$ is in $T_q$, thus the composed map 
$X_\lambda \to T_{q+1}$ admits by construction 
a factorization through an acyclic $F$. We now 
build the pushout $G=F \coprod_{X_\lambda} Y$, 
that is acyclic : $X_\lambda \to Y' \to Y$ is a 
homology equivalence, thus $F \to G$ is one also. 
By construction, the map $f$ factorizes through $G$,
hence the colimit-object $Z$ is $\A $-local. 
Finally, consider $0 \to \build
{\oplus}_0^\infty Z_n/X \to \build{\oplus}_0^\infty
Z_n/X \to Z/X \to 0$ this short exact sequence 
shows that the colimit-map $X \to Z$ is a 
$\tilde\A $-equivalence. $\blacksquare$  
\end{demo}

\begin{lemme}.\label{quatre}\\
Let $X$ be an $\A $-local $\cD$-complex. 
Then $X$ is also $\tilde \A $-local.
\end{lemme}
\begin{demo}\\
\vspace*{-.6cm}\\
Let ${\call C}$ be the class of all $\cD$-complexes $Y$,
such that, for all $L$ in $\L\A$, every map $Y 
\to L$ factorizes through an acyclic. This class 
contains $\A $. Like $\L\A$, it's stable under 
desuspension. It's stable under direct sum : actually, 
let $Y=\oplus Y_i$ and $f=\oplus f_i:Y_i \to L$ 
with $L$ in $\L\A$; each $f_i$ factorizes through an 
acyclic $E_i$, so $f$ factorizes through $E=\oplus E_i$ 
acyclic. Let's prove that ${\call C}$ is stable under the 
Coker of a cofibration : let $0 \to X \to Y \to Z \to 0$ 
be a short exact sequence, with $X$ and $Y$ in ${\call C}$; 
let $L$ be in $\L\A$ and consider $f:Z \to L$. The 
composed map $Y \to Z \to L$ factorizes through $E$ acyclic.
We can choose $E$ such that $E \to L$ is surjective, and 
denote $L'$ the Kernel. As $\L\A$ is an exact class, 
$L'$ is in $\L\A$, thus the map at the level of the 
Kernels $X \to L'$ factorizes through $F$ acyclic. Let $E'$ 
be the direct sum of $E$ and an acyclic containing $F$ by 
a cofibration. Denoting $G=Coker[F \to E']$, that is also 
acyclic, we can complete the diagram of compatible short 
exact sequences : 
$$\xymatrix@-1pc{
0 \ar[r] & L' \ar[r] & E \ar[r] & L \ar[r] & 0 \\
0 \ar[r] & F \ar[r] \ar[u] & E' \ar[r] \ar[u] &
G \ar[r] \ar[u] & 0 \\
0 \ar[r] & X \ar[u] \ar[r] & Y \ar[r] \ar[u] &
Z \ar[r] \ar[u] & 0}$$ 
\mbox{The class ${\call C}$ is exact, stable under direct
sum, and contains $\A $. Thus ${\call C}$ contains 
$\tilde \A $ too. $\blacksquare$} 
\end{demo}

\begin{lemme}.\\
Let $f:X \to Y$ be a morphism of $\cD$-complexes,
with $X \in \B $ and $Y \in \tilde \A $.
\\Then $f$ admits a factorization through a 
$\cD$-complex $Z \in \A $.
\end{lemme}
\begin{demo}\\
\vspace*{-.6cm} \\
Let ${\call C}$ be the class of all $\cD$-complexes $Y$ in 
$\tilde \B $ such that, for all $X$ in $\B $, 
every map $X \to Y$ factorizes through an object in 
$\A $. This class contains $\A $. Like $\tilde 
\B $, it's stable under desuspension. Let's prove 
that ${\call C}$ is stable under direct sum : let $Y_i$ be 
objects in ${\call C}$ and consider a map $f=\oplus f_i :
X \to Y=\oplus Y_i$. As $\B $ is h-compact, there 
exists a homology equivalence $X' \build{\to}_{}^{\simeq} X$
such that $X' \to X \to Y$ factorizes through a finite sum 
$Y'=\build{\oplus}_{0}^{p} Y_i$. The map $X' \to Y'$ 
factorizes through an object $Z$ in $\A $, because 
the class ${\call C}$ is obviously stable under finite 
direct sum ( in the finite case, we have ``product=sum'' ). 
The map $f$ factorizes through the pushout $U=X\coprod_{X'} 
Z$ which is also in $\A $ : actually, $X' \build{\to}
_{}^{\simeq} X$ is a homology equivalence, so $Z \build{\to}
_{}^{\simeq} U$ is one too; then $Z$ is in $\A $ 
implies that $U$ is also ! \mbox{We can sum up the situation
in the following diagram :} 
$$\xymatrix@-1.3pc{
X' \ar[r] \ar[d]_\simeq & 
(Z \in \A ) \ar@{-->}[d]_\simeq \ar@/^/[dr] \\
X \ar@{-->}[r] \ar@/_/_f[rrd] & 
U \ar@{..>}[dr] & (Y'=\build{\oplus}_0^p Y_i) \ar[d] \\
&& (Y=\oplus Y_i)}$$
Let's now prove that the class ${\call C}$ is stable under 
the Coker of a cofibration. Let $0 \to X \to Y \to Z \to 0$
be a short exact sequence, with $X$ and $Y$ in ${\call C}$.
Let $U$ be in $\B $ and consider a map $f:U \to Z$.
Let $E$ be an acyclic that maps onto the pullback $Y \prod_Z
U$. Denote $K=Ker[E \to U]$. As $X$ is in ${\call C}$, the 
map $K \to X$ factorizes through an object $X_1$ in 
$\A $. As $Y$ is in ${\call C}$, the map from the 
pushout $X_1 \coprod_K E$ to $Y$ factorizes through an 
object $Y_1$ in $\A $. We can suppose that the map 
$X_1 \to Y_1$ is a cofibration, and denote $Z_1$ the 
Cokernel ( if need be, we can add acyclics ). By 
factorization of the Cokernels, the map $U \to Z$ factorizes
through $Z_1$, as we can see on the following diagram : 
$$\xymatrix@-1pc{
0 \ar[r] & X \ar[r] & Y \ar[r] & Z \ar[r] & 0 \\
0 \ar[r] & X_1 \ar[r] \ar[u] & Y_1 \ar[r] \ar[u] &
Z_1 \ar[r] \ar[u] & 0 \\
0 \ar[r] & K \ar[u] \ar[r] & E \ar[r] \ar[u] &
U \ar[r] \ar[u] & 0}$$    
\mbox{The class ${\call C}$ is exact, stable under direct 
sum, and contains $\A $; thus ${\call C}$ contains 
$\tilde \A $. \hfill $\blacksquare$} 
\end{demo}

\begin{lemme}.\\
Let $\B '$ be the class of $\cD$-complexes 
$X$, such that there exists a $\cD$-complex $U \in 
\B $ and a $\tilde \A $-equivalence 
$U \to X$. Then the canonical inclusion $\B  
\subset \B '$ induces a homotopy equivalence 
of spectra : $K(\B,\A)\simeq K(\B',\tilde \A)$.
\end{lemme}
\begin{demo}
We verify the hypothesis of Waldhausen's approximation 
theorem for the inclusion functor $F: \B  \subset 
\B '$. [App1] : Let $X$ be an object in $\B $
such that $F(X)$ is in $\tilde \A $. The identity 
map $Id:X \to X$ factorizes through an object in 
$\A $. Thus there exists $X'$ in $\B $ such 
that $X \oplus X'$ is in $\A $. As $\A $ is 
stable under direct summand, $X$ is already in 
$\A $. [App2] : Let $X$ be an object in $\B $
and $f:X \to Y$ any map in $\B '$. There exists $Z$
in $\B $ and an $\tilde \A $-equivalence 
$Z \to Y$. Let $E$ be an acyclic that contains $Z$ by 
a cofibration, and denote $U=Coker[Z \to Y \oplus E]$. 
The $\cD$-complex $U$ is in $\tilde \A $, and the 
map $X \to Y \to U$ factorizes through an object $V$ in 
$\A $. We build the pullback $X'=V \coprod_U 
(Y \oplus E)$ : the map $V \to U$ is an $\tilde 
\A $-equivalence, thus $X' \to Y \oplus E$ is one 
also; finally, $X'$ is the extension of $Z$ by $V$ ( two 
objects in $\B $ ), thus $X'$ is in $\B $ 
also. The following diagram proves the surjectivity 
hypothesis, and then Waldhausen's approximation theorem 
gives the equivalence of K-theory spectra : \hfill 
$K(\B,\A) \simeq K(\B',\tilde \A)$. 
$$\xymatrix@-1pc{
& Z \ar@{=}[r] \ar@{ >->}[d] & Z \ar@{ >->}[d] \\
X  \ar[r] \ar[dr] & X' \ar[r] \ar@{->>}[d] &
Y \oplus E \ar@{->>}[d] \\
& V \ar[r] & U}$$
Remark that the class $\B '$ contains the acyclics.
It's stable under desuspension ( because $\B $ and
$\tilde \A $ are ). Finally, it's stable under the 
Coker of a cofibration : let $0 \to X_1 \to X_2 \to X_3
\to 0$ be a short exact sequence, with $X_1$ and $X_2$ in
$\B '$. There exists an $\tilde 
\A $-equivalence $U_1 \to X_1$ with $U_1$ in 
$\B $. We apply the axiom [App2] to the map $U_1 \to
X_1 \to X_2$ to find a map $u:U_1 \to U_2$ in $\B $ 
and an $\tilde \A $-equivalence $U_2 \to X_2$ such 
that the square commutes. Denote $U_3=Coker^h(u)$ in 
$\B $, we obtain by passing to the Cokernel the 
$\tilde \A $-equivalence $U_3 \to X_3$ wanted : 
actually, $\B $ and $\tilde \A $ are stable 
under '2/3' axiom. Hence the class $\B '$ is exact.
\hfill $\blacksquare$ 
\end{demo}

\begin{lemme}.\\
Let ${\call E}$ be the category of short exact 
sequences of $\cD$-complexes $\xymatrix@-.8pc{
X \ar@{ >->}[r] & Y \ar@{->>}[r] & Z}$ with 
$X \in \tilde \A $ and $Z \in \L$. 
If we choose the homology equivalences at the level 
of the quotient term $Z$ to be the class 
$w{\call E}$ of weak equivalences, then ${\call E}$ 
has a structure of Waldhausen category. Let $F$ be 
the exact functor that sends an exact sequence 
$\xymatrix@-.8pc{X \ar@{ >->}[r] & Y \ar@{->>}[r] & 
Z}$ to its middle term $Y$. By Waldhausen's 
approximation theorem, the functor $F$ induces a 
homotopy equivalence of spectra :
$K({\call E}) \simeq K(\B',\tilde \A)$. 
\end{lemme}
\begin{demo}
Let $F$ be the functor that sends each short exact 
sequence to its middle term. We shall verify the 
hypothesis of Waldhausen's approximation theorem 
for $F$. [App1] : Let $\xymatrix@-1pc{X \ar@{ >->}[r] &
Y \ar@{->>}[r] & Z}$ be a short exact sequence in 
${\call E}$. If $Y$ is in $\tilde \A $, then $Z$ is 
in $\L$ and in $\tilde \A $. By lemma 
\ref{quatre}, $Z$ is $\tilde \A $-local, and the 
identity $Id:Z \to Z$ factorizes through an acyclic. Thus 
$Z$ is acyclic. We apply this argument to $Coker^h(f)$
for any map $f$ in ${\call E}$ to deduce : $F(f)$ is 
in $\tilde \A $ if and only if $f$ is a weak 
equivalence in ${\call E}$. [App2] : Let $\xymatrix@-1pc{(X 
\ar@{ >->}[r] & Y \ar@{->>}[r] & Z)}=[S]$ be a short exact 
sequence in ${\call E}$, and $f:Y \to Y_0$ any map in 
$\B '$. We build the pushout $Z_0=Z \coprod_Y Y_0$. 
As $\B $ is h-small and h-compact, so is $\A $; 
now by lemma \ref{trois}, there exists an $\tilde 
\A $-equivalence $Z_0 \to Z'$ with $Z'$ an 
$\A $-local object. We add an acyclic $E$ that surjects
onto $Z'$ to obtain a surjective map $\xymatrix@-1pc{
(Y'=Y_0 \oplus E) \ar@{->>}[r] &Z'}$, denote $X'$ its Kernel,
we get a new short exact sequence in ${\call E}$ : 
$\xymatrix@-1pc{(X' \ar@{ >->}[r] & 
Y' \ar@{->>}[r] & Z')}=[S']$, doted with a map of 
short exact sequences $\Phi:S \to S'$, and the 
\mbox{map $f$ factorizes as $F(\Phi)$ followed by a homology 
equivalence (cf the following diagram). $\blacksquare$} 
$$\xymatrix@-1.2pc{
[S] \ar[dd]_\Phi & (X \in \tilde\A ) \ar@{ >->}[r] 
\ar@{..>}[dd] & (Y \in \B ') \ar@{->>}[r] \ar[d]_f &
(Z \in \L) \ar@{..>}[d] \\
&& (Y_0 \in \B ') \ar@{..>}[r] \ar[d]_{(Id \oplus 0)}
^\simeq & Z_0 \ar[d]^{\tilde \A -eq} \\
[S'] & X' \ar@{ >->}[r] & (Y'=Y_0 \oplus E) \ar@{->>}[r] 
\ar[d]_{(Id \oplus 0)}^\simeq & Z' \\
&& Y_0}$$     
\end{demo}

\begin{lemme}.\\
Pose $\L=\B' \cap \L\A$. Then $\L$ is an 
exact class of $\cD$-complexes. Let $G$ be the exact
functor that sends an exact sequence 
$\xymatrix@-.8pc{
X \ar@{ >->}[r] & Y \ar@{->>}[r] & Z}$ 
to its quotient term $Z$. By Waldhausen's additivity
theorem, the functor $G$ induces a homotopy 
equivalence of spectra : 
$K({\call E}) \simeq K(\L)$.
\end{lemme}
\begin{demo}
The class $\L$ is the intersection $\B' \cap \L\A$, 
hence $\L$ is exact. By the additivity 
theorem, the K-theory spectrum $K({\call E})$
\mbox{is equivalent to the product $K(\tilde \A , 
\tilde \A ) \times K(\L)$ and thus to 
$K(\L)$. \hfill $\blacksquare$} \\
\end{demo}

\no
This ends the proof of the Main Theorem. \hfill 
$\blacksquare$

\section{Excision Theorem for Transverse Classes}

\begin{defin}.\\
Let $\A,\B$ be two exact classes of 
$\cD$-complexes. \\We shall say that {\bf 
``$\A $ is transverse to $\B $''}, 
denoted by $\A \pitchfork \B$, if every map 
$f:A \to B$ with $A \in \A $ and 
$B \in \B $ admits a factorization through a 
common term $C \in \A \cap \B$. In that case, 
we shall note $\A +\B$ the smallest exact 
class of $\cD$-complexes containing $\A $ 
and $\B $.
\end{defin}

\begin{prop}.\\
Let $\A,\B$ be two exact classes of 
$\cD$-complexes. Suppose that $\A \pitchfork \B$. Let $\C$ 
be the class of $\cD$-complexes $X$ ( non-necessarily h-finite ),
such that there exists $A \in \A $, $B \in \B $, 
and a short exact sequence :
$\xymatrix@-1pc{ B \ar@{ >->}[r] & A \ar@{->>}[r] &
X}$. This class $\C$ has a natural 
cylinder-functor, and verifies the '2/3 axiom'. 
As $\C$ evidently contains $\A $ and 
$\B $, is exact, and is the smallest for these
properties, thus we have : $\C=\A+\B$. 
\end{prop} 
\begin{demo}
$\bullet$ With the definition above, we show the natural 
inclusions $\A \subset \A+\B$ by the short exact 
sequence $0 \to A \to A$; and $\B \subset \A+\B$ by 
the short exact sequence $s^{-1}B \to 0 \to B$. The first 
reflex is to verify if a short exact sequence `in the other
way' $A \to B \to X$ gives also an object in $\A+\B$. 
To this aim, we factorize $A \to B$ by $C$ in $\A \cap \B$, 
and write the diagram :   
$$\xymatrix@-1pc{
A \ar[r] \ar[d] & C \ar[r] \ar[d] & \A  \ar[d] \\
B \ar@{=}[r] \ar[d] & B \ar[r] \ar[d] & 0 \ar[d] \\
X \ar[r] & \B  \ar[r] & \A }$$
$\bullet$ Let's look now at the four elementary cases of 
extension : we shall take the following notations \\
\mbox{$(X \in \A+\B), \, (A,A',A'' \in \A ), \, 
(B,B',B'' \in \B )$, and $(C \in \A \cap \B)$ 
obtained by factorization.} 
$$\xymatrix@-1pc{
B \ar[r] \ar@{=}[d] & A'' \ar[r] \ar[d] & ? \ar[d] &
0 \ar[r] \ar[d] & B' \ar@{=}[r] \ar[d] & B' \ar[d] &
\A  \ar[r] \ar[d] & A' \ar[r] \ar[d] & C \ar[d] &
\A  \ar[r] \ar[d] & ? \ar[r] \ar[d] &\B  \ar[d]\\
B \ar[r] \ar[d] & A \ar[r] \ar[d] & X \ar[d] &
A \ar[r] \ar@{=}[d] & X \ar[r] \ar[d] & B \ar[d] &
A \ar[r] \ar[d] & X \ar[r] \ar[d] & B \ar[d] &
A \ar[r] \ar[d] & X \ar[r] \ar[d] & B \ar[d]\\
0 \ar[r] & A' \ar@{=}[r] & A' &
A \ar[r] & ? \ar[r] & B'' &
\A  \ar[r] & ? \ar[r] & \B  &
C \ar[r] & B' \ar[r] & \B }$$
These four diagram must be read vertically, and show that 
the class $\A+\B$ is stable under : the Kernel of a map
going to $\A $; the Cokernel of a map coming from 
$\B $; the Cokernel of a map coming from $\A $;
finally, the Kernel of a map going to $\B $. By 
composing these four elementary operations, we shall treat 
the '2/3' axiom : consider $(A \to X \to B),(A' \to Y \to B')
\in \A+\B$.
$$\xymatrix@-1.8pc{
A \ar[rr] \ar@{=}[dd] && X \ar[rr] \ar[dd] && B \ar[dd] \\
&&& \textrm{\textcircled{p}} \\
A \ar[rr] \ar[dd] && Y \ar[rr] \ar[dd] && Z \ar[dd] \\
\\
0 \ar[rr] && K \ar@{=}[rr] && K}  \hspace*{3cm} 
\raisebox{-0.95cm}{\parbox{5.4cm}{
We consider the pushout $Z$ here. The third case above 
shows that $A \in \A $ and $Y \in \A+\B$ 
imply $Z \in \A+\B$. Then the second case above 
shows that $Z \in \A+\B$ and $B \in \B $ 
imply $K \in \A+\B$. Thus the '2/3' axiom holds.}}$$
$\bullet$ For the cylinder-functor $T(f:X \to Y)$, we shall 
consider the short exact sequence \\$\xymatrix@-1pc{(Y 
\ar@{ >->}[r] & T(f) \ar@{->>}[r] & sX)}$; and then apply 
the stability under suspension and '2/3' axiom. \hfill 
$\blacksquare$     
\end{demo}

\begin{theo} {\em [Excision] [Vogel,Bihler]} \\
Let $\A,\B$ be two exact h-small classes of 
$\cD$-complexes. Suppose that $\A \pitchfork \B$. \\
Then we have the homotopy equivalence of spectra : 
$K(\A+\B,\A) \simeq K(\B,\A \cap \B)$.
\end{theo}
\begin{demo}
We apply Waldhausen's approximation theorem to the 
inclusion $(\B,\A \cap \B) \subset (\A+\B,\A)$. 
We must therefore define the weak equivalences  
\mbox{$(f:X \to Y) \in v(\A+\B) \Leftrightarrow 
Coker^h(f) \in \A $}. Similarly, $(g:B \to B') \in 
w\B  \Leftrightarrow Coker^h(g) \in \A \cap \B$.
The only non-trivial hypothesis that we must verify is the 
surjectivity axiom [App2]. Consider on the left-side 
$B' \in \B $ and on the right-side 
$\xymatrix@-1pc{(A \ar@{ >->}[r] & X \ar@{->>}[r] & B)} 
\in \A+\B$ doted with a map $f:X \to B'$.
We want to factorize $f$ through a map in $\B $, modulo 
a $v$-equivalence. The composed map $A \to X \to B'$ factorizes 
through an object $C \in \A \cap \B$. We build the pushout
$Y=X \coprod_A C$. The map $X \to Y$ is a $v$-equivalence :
actually, $Coker^h(X \to Y)=Coker^h(A\to C) \in \A $. 
Moreover, the map $f$ factorizes as $X \to Y \to B'$. Finally,
the object $Y$ is in $\B $ ( it's the extension of $C$
by $B$, both in $\B $, which is stable under '2/3' ). 
The axiom [App2] is hence verified, and Waldhausen's 
approximation theorem gives : $K(\B, \A \cap \B) \simeq 
K(\A+\B,\A)$. This situation is summed up in the following 
diagram : \hfill $\blacksquare$
$$\xymatrix@-1.7pc{
A \ar@{ >->}[rr] \ar[dd]_v && X \ar@{->>}[rr] 
\ar@{-->}[dd]_v \ar@/^/[dddr] && B \ar@{=}[dd] \\ 
& \lrcorner \\ C \ar@{ >-->}[rr] \ar@/_/[rrrd] && 
Y \ar@{..>}[rr] \ar@{..>}[dr] && B \\ &&& B'}$$
\end{demo}

\no
{\bf Remark :} \\
Even if the transversality hypothesis is not symmetric, 
the conclusion of our excision theorem is ! Actually, 
consider the following diagram where the lines \& columns 
are homotopy fibrations :
$$\xymatrix@-1.8pc{
K(\B,\A \cap \B) \ar@{ >->}[rr] && K(\A+\B,\A) \\ \\
K(\B ) \ar@{->>}[uu] \ar@{ >->}[rr] && K(\A+\B) 
\ar@{->>}[uu] \ar@{->>}[rr] && K(\A+\B,\B) \\
& \textrm{\textcircled{p}} \\
K(\A \cap \B) \ar@{ >->}[uu] \ar@{ >->}[rr] &&
K(\A ) \ar@{ >->}[uu] \ar@{->>}[rr] && 
K(\A,\A \cap \B) \ar@{ >->}[uu]}$$
If one of the two inclusions ( above or on the right ) gives 
an equivalence in K-theory, then the square \textcircled{p} 
is 'exact' ( ie Cartesian \& Cocartesian ), and this implies 
that the other map gives a homotopy equivalence between the 
K-theory spectra. \mbox{Anyway, we can deduce a Mayer-Vietoris 
long exact sequence $(\forall i \in \Z)$:} 
$$\ldots \to K_i(\A \cap \B) \to K_i(\A ) \oplus 
K_i(\B ) \to K_i(\A+\B) \to 
K_{i-1}(\A \cap \B) \ldots$$
{\bf Example :} \\
Let $R$ be a noetherian regular ring. For any ideal, we have a 
notion of order defined via the Krull-dimension. For each 
complex $C_*$, consider the annulator ideal ${\call I}=
\{ a \in R \; | \; \xymatrix@-1pc{C_* \ar[r]^{.a} & C_*} \sim 0 
\}$ : we get a notion of order for complexes. Now take two 
prime integers $p \wedge q=1$; define \mbox{$\A = 
\{ C_* \; | \; \exists m \ge 1, order( H_*(C))=p^m \}$} and 
$\B = \{ C_* \; | \; \exists n \ge 1, 
order( H_*(C))=q^n \}$. Then we can verify that 
$\A \pitchfork \B$, and thus we can apply the excision 
theorem, knowing that $\A \cap \B$ contains only the 
acyclics : thus the equivalence $K_i(\A+\B) \simeq 
K_i(\A ) \oplus K_i(\B), \forall i \in \Z$. 
\section{Application to the $\tilde K{\cN}il$ groups}

We shall introduce the structures we need to study Waldhausen's 
$\tilde K{\cN}il$ groups in the setting of Vogel's [Conjecture].
Hereafter, $R$ will be a fixed regular ring in the sense of 
Vogel and $S$ a fixed $R$-bimodule, flat on the left ( cf 
\cite{bil1} for precisions ). After a few technical definitions,
we shall apply our localization and excision theorems on the 
class ${\call B}_n$ and ${\call A}_n$ of nilpotent 
$\cD_0$-complexes, and gives examples of concrete calculus of 
local objects.     
\begin{defin}.\\
Let $\cD_0$ be the diagram of bimodules given by :\\
$\bullet$ The oriented graph $\Gamma$ is the set $\N$ of 
natural integers, doted with maps \\ \hspace*{.2cm }
$(\lambda : n \to n+1)$ and $(\alpha : n+1 \to n)$ such 
that $\alpha \circ \lambda=\lambda \circ \alpha$. \\
$\bullet$ To each vertex $n \in \N$ we associate the ring $R$. 
\\\mbox{$\bullet$ To each map $\lambda$ we associate the 
bimodule $R$, to each map $\alpha$ the bimodule $S$.}
\end{defin}
\begin{defin}. \\ \vspace*{-.7cm} \\
 
$\left \llbracket 
\parbox[l]{12cm}{
For all $n \ge 0$, let ${\call B}_n$ be the class of 
$\cD_0$-complexes \\
$\xymatrix@-1pc{
B_*= &(0 \ar[r] & B_1 \ar[r]_{\lambda_1} 
\ar@/_1pc/@{..>}[l]_{\alpha_1}
&B_2 \ar[r]_{\lambda_2} \ar@/_1pc/@{..>}[l]_{\alpha_2} &
B_3 \ar@/_1pc/@{..>}[l]_{\alpha_3} & \ldots & 
B_n \ar@{=}[r]_{\lambda_n} \ar@/_1pc/@{..>}[l]_{\alpha_n} &
B_n \ar@{=}[r]_{\lambda_{n+1}} 
\ar@/_1pc/@{..>}[l]_{\alpha_{n+1}} & B_n & \ldots)}$ \\
where $B_0=0$ and each $B_i$ is a h-finite complex \\( of 
projective right $R$-modules ) \\
each $\lambda_i:B_i \to B_{i+1}$ is a chain map, and a 
cofibration \\( We suppose furthermore that $\lambda_i$ 
is a weak equivalence for all $i \ge n$ ) \\
each $\alpha_i:B_i \to B_{i-1} \otimes _R S$ is a chain map \\
( If all $\alpha_i$ are surjective, the $\cD_0$-complex $B_*$
will be called 'reduced' ) \\
finally, we ask that the two composed maps $B_i \to B_{i+1} 
\to B_i \otimes S$ \\and $B_i \to B_{i-1} \otimes S \to B_i 
\otimes S$ verify $\alpha \circ \lambda=\lambda \circ \alpha$.}
\right .$ 
\end{defin}
\begin{defin}.\\
$\bullet$ For all $n \ge 0$, let ${\call A}_n$ be the class 
of $B_* \in {\call B}_n$ such that $B_n$ is contractible. \\
$\bullet$ The inclusion functors ${\call B}_n 
\subset {\call B}_{n+1}$ for all $n \ge 0$ define 
${\call B}= \varinjlim {\call B}_n$. \\
$\bullet$ The inclusion functors ${\call A}_n 
\subset {\call A}_{n+1}$ for all $n \ge 0$ define 
${\call A}=\varinjlim {\call A}_n$.   
\end{defin}

\begin{theo}.\\
Let $R$ be a right regular ring (in the sense of Vogel ), 
and let $S$ be any $R$-bimodule, flat on the left. 
Then we have the following homotopy equivalence
of non-connective K-theory spectra : \hfill 
$K{\cN}il(R;S) \simeq K({\call B,A}) \simeq \varinjlim 
K({\call B}_n,{\call A}_n)$.
\end{theo}

\no
We refer the interested reader to \cite{bil2} for a complete 
proof. \\

\no
{\bf Calculus of local objects :} \\
\ding{182} We shall first treat the easy graded case ${\call B}_n 
\subset {\call B}_{n+1}$. 
\begin{prop}.\\
Let $C_*$ be a reduced $\cD_0$-complex in $\tilde {\call B}$. \\
Then ``$C_*$ is ${\call B}_n$-local'' if and 
only if ``$C_i$ is contractible $\forall 0 \le i \le n$''.
\end{prop}
\begin{demo}
$\bullet \; \underline{(i) \to (ii)}$ : We consider the test-objects 
$D_*=g_m(D)=(0 \to 0 \to \ldots 0 \to D = D = \ldots)$ beginning with
$m$ zeros, where the complex $D=(\ldots \to 0 \to R \to 0 
\to \ldots)$ is the base ring concentrated in degree $0$, 
and every index $1\le m \le n$. Let $\cH om(D_*,C_*)$ denote
the set of maps $f_*:D_* \to C_*$ respecting the graduations 
$f_k:D_k \to C_k$, the structural maps $(\lambda,\alpha)$ of 
$\cD_0$-complexes, but a priori not the degrees, and not the 
differentials ( we call these maps 'algebraic morphisms' ). 
It's a graded differential $R$-module : $\cH om(D_*,C_*)_p
=\{f:D_k^l \to C_k^{l-p} \}$ and following Leibniz rule : 
$(df)=d \circ f - (-1)^{deg(f)} f \circ d$. By $R$-linearity, 
a morphism $f_*:D_* \to C_*$ is given by the image $x=f(\bI _R) 
\in C_m$. To respect the $\lambda$ maps, we pose $f_k=0$ for 
$k<m$ and $f_{m+k+1}=\lambda_{m+k} \circ \ldots \lambda_m 
\circ f$. To respect the $\alpha$ maps, we have the condition 
$\alpha_m(x)=0$. Then we write $\alpha_{m+k} \circ f_{m+k} 
(\bI _R)= \lambda_{m+k-2} \circ \ldots \lambda_m \circ \lambda_
{m-1} \circ \alpha_m(x)=0$ due to the relation $\alpha \circ 
\lambda = \lambda \circ \alpha$ available in any $\cD_0$-complex.
Thus we have an isomorphism of graded $R$-modules : $\cH om
(D_*,C_*) \simeq Ker(\alpha_m)$ that sends each $f$ on $x=f(\bI 
_R)$. Furthermore, $df=0$ corresponds to a cycle $dx=0$, and 
$f=dg$ can be interpreted as : $x=f(\bI _R),y=g(\bI _R) \in 
Ker(\alpha_m)$ and $x=dy$ [ ie it's a boundary ]. Hence it's 
an isomorphism of differential graded $R$-modules. Let now 
$f$ be in $\cH om(D_*,C_*)$ such that $df=0$ : it's a chain
morphism of degree $p$. As $C_*$ is ${\call B}_n$-local and 
$D_* \in {\call B}_n$, the map $f$ factorizes through an acyclic 
$E_*$. Thus $f=dg$ is a boundary. Conclusion : the complex 
$Ker(\alpha_m)$ is acyclic. A finite induction on $i$ shows 
that each complex $C_i$ is acyclic for $i \le n$ : actually, 
it's true for $i=0$; now look at the short exact sequence 
$0 \to Ker(\alpha_i) \to C_i^* \to C_{i-1}^* \otimes S \to 0$ 
( here the object $C_*$ is reduced ). Then $C_{i-1}^*$ is acyclic
by induction hypothesis, thus $C_{i-1}^* \otimes S$ is acyclic 
too, because $S$ is flat on the left; $Ker(\alpha_i)$ is acyclic
due to the condition imposed by the test-object $g_i(D)$; 
hence $C_i^*$ is acyclic also. We have proven condition $(ii)$.\\
$\bullet \; \underline{(ii)\to(i)}$ : Conversely, we shall now verify
that every object of this type is ${\call B}_n$-local. Let $f:D_*\to 
C_*$ be any map from any $D_* \in {\call B}_n$ to $C_*$ with $C_i$ 
contractible $\forall 0 \le i \le n$. We shall build by induction 
on $i$ a contractible object $E_*$ through which $f$ admits a 
factorization. For $i \le n$, we pose $E_i=C_i$ doted with the same 
$(\lambda,\alpha)$ structural maps. For $i \ge n$, we pose 
$E_{i+1}=E_i \coprod_{D_i} D_{i+1}$. So now we must define the maps 
$\lambda_i: E_i \to E_{i+1}$ and $\alpha_{i+1}:E_{i+1} \to E_i 
\otimes S$ in a compatible way with the maps on $D_*$ and $C_*$. 
As $D_i \to D_{i+1}$ is a cofibration and a homology equivalence 
( because $D_*$ is in ${\call B}_n$ ), the pushout gives $\lambda_i:
E_i \to E_{i+1}$ cofibration and homology equivalence, compatible 
with $D_*$. Next consider the following pushout, that will define the
map $E_{i+1} \to C_{i+1}$ ( by the universal property ) compatible
with the maps $\lambda$ :            
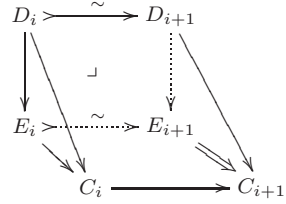
$$\xymatrix@-1.2pc{
D_i \ar@{ >->}[rr]^\sim \ar[dd] \ar[rddd] && D_{i+1} \ar@{..>}[dd]
\ar[rddd] \\ & \lrcorner \\
E_i \ar@{ >..>}[rr]^\sim \ar[rd] && E_{i+1} \ar@{=>}[rd] \\
& C_i \ar[rr] && C_{i+1}}$$   
Actually, the exterior ``square'' commutes, because the initial map
$f$ commutes with $\lambda$. We must now define $(\alpha_{i+1}:
E_{i+1} \to E_i \otimes S)$ : on $D_{i+1}$, take  
$D_{i+1} \to D_i \otimes S \to E_i \otimes S$; on $E_i$, take 
the composition $E_i \to E_{i-1} \otimes S \to E_i \otimes S$. 
\mbox{These maps coincide on $D_i$ as suggested below :}
$$[\xymatrix@-1pc{
{.} \ar@<-2pt>_1[r] \ar[d]_3 & {.} \ar@<-2pt>_2[l] & = & 
{.} \ar@<-2pt>_2[r] & {.} \ar@<-2pt>_1[l] \ar^3[d] & = &
{.} \ar_2[d] & {.} \ar_1[l] & = &
& {.} \ar^1[d] & = & 
{.} \ar_1[d] 
\\{.} &&&&{.} &&{.} \ar_3[r]& {.} && {.} \ar@<-2pt>_3[r] & {.} 
\ar@<-2pt>_2[l] && {.} \ar@<-2pt>_2[r] & {.} \ar@<-2pt>_3[l]}]$$ 
Thus we can define $\alpha_{i+1}$ on $E_{i+1}$ by the universal 
property of the pushout, and it will be compatible with that 
defined on $D_*$. It remains to verify the compatibility with $C_*$
: on $D_{i+1}$, it's just the compatibility of the initial map $f$
with $\alpha$; on $E_i$, we go back to the induction hypothesis, 
with exactly the same manipulation as above, but between $E_i$ and
$C_{i+1}$. Hence we have built by induction on $i$ a contractible 
$\cD_0$-complex $E_*$ through which $f:D_* \to C_*$ admits a 
factorization. The objects verifying $(ii)$ are thus well ${\call
B}_n$-local. \hfill $\blacksquare$ \\
\end{demo}

\no
We can now identify the local objects, via Waldhausen's approximation
theorem, to obtain a substitute of Quillen's d\'evissage for the 
class ${\call B}$ of $\cD_0$-complexes. 

\begin{prop}.\\
Let $\L$ be the class of reduced ${\call B}_n$-local 
$\cD_0$-complexes $L$, such that there exists $B$ in 
${\call B}_{n+1}$ and a $\tilde {\call B}_n$-equivalence from 
$B$ to $L$. By Waldhausen's approximation theorem, the exact 
projection functor $F:\L \to \C_R$ defined by 
$L_* \mapsto L_{n+1}$ induces a homotopy equivalence of 
K-theory spectra : $K({\call B}_{n+1},{\call B}_n) \simeq K(R)$.  
\end{prop}
\begin{demo}
$\bullet$ By use of the main theorem, we know the equivalence :
$K({\call B}_{n+1}^{red}, {\call B}_n^{red}) \simeq K(\L)$. 
It now remains to identify the objects $L$ in $\L$ : by the 
proposition above, $L_i$ is contractible $\forall 0 \le i \le n$. 
Moreover, we have a short exact sequence $\xymatrix@-1pc{(A 
\ar@{ >->}[r] & B \ar@{->>}[r] & L)}$ with $A$ in 
$\tilde {\call B}_n$ and $B$ in ${\call B}_{n+1}$. In particular, 
the $\lambda_i$ are homology equivalence $\forall i>n$. Let's prove 
that $L_{n+1} \in \C_R$, and thus that the projection functor $F$ 
is well defined : consider the exact sequence $(A_n \to B_n \to *)$; 
$B_n$ is in $\C$, thus $A_n$ too. Now look at the exact sequence 
\mbox{$(A_{n+1} \to B_{n+1} \to L_{n+1})$; as $A_{n+1}=A_n$ is in 
$\C_R$ and $B_{n+1}$ too, hence $L_{n+1}$ is in $\C_R$ also !} \\   
$\bullet$ The projection functor $F$ sends a homology equivalence 
in the class of local objects to a homology equivalence in $\C_R$.
Conversely, let $f:A_* \to B_*$ between local objects, such that 
$f_{n+1}$ be a homology equivalence. As $A_i \simeq A_{i+1}$ and 
$B_i \simeq B_{i+1}$ for all index $i>n$ ( these are homotopy 
equivalences because we work in the class $\C_R$ ), then 
$f_i \simeq f_{n+1}$ is a homology equivalence. For the index 
$i \le n$, all complexes $A_i$ and $B_i$ are contractible, 
so $f_i$ is a homology equivalence ! Thus $f_*=(f_i)$ is a 
homology equivalence, and the axiom [App1] is verified. \\
$\bullet$ Let's now study the surjectivity [App2] : let $C_*$ be 
a local object described in the proposition above, and let $f:C_{n+1}
\to D$ be any map in $\C_R$. We want to build a $\cD_0$-complex $D_*$
local and a map $f_*$ that lift $f$. By a factorization through the 
cylinder-functor on $\C_R$, we can suppose that $f$ is a cofibration.
We shall build $D_*$ by induction on $i$. For the index $i<n$, pose 
$D_i=C_i$ with the same structural maps $(\lambda,\alpha)$. For 
$i=n$, pose $D_n=C_n$ with the same $\alpha_n$, but here $\lambda_n
=f \circ \lambda_n$. Problem : can we define $\alpha_{n+1}:D \to C_n
\otimes S$ that lift the map $\alpha_{n+1}$ already defined on 
$C_{n+1}$ ? As $C_n$ is acyclic, and $S$ flat on the left, then 
$C_n \otimes S$ is acyclic also : there is no obstruction to lift 
along the cofibration $f$. Then it remains all the right-side of 
$D_*$ to build. For this, we use the same argument as in the proof 
$(ii) \to (i)$ above : the pushout ! For $i >n$, pose $D_{i+1}=D_i 
\coprod_{C_i} C_{i+1}$, and now all structural maps $(\lambda,
\alpha)$ are defined functorially ! So, we have built a map of 
$\cD_0$-complexes $f_*$ that lifts $f$ ( followed eventually by 
a homotopy equivalence after the factorization through the cylinder 
\ldots ) and the multicomplex $D_*$ we've constructed has the same 
properties as $C_*$ : contractible for the index $i \le n$, then 
the $\lambda_i$ cofibrations and homology equivalences make their 
way in the pushout. So $D_*$ is local. Moreover $D_{n+1}=D_n \in 
\C_R$, thus the constructed object $D_*$ is in $\L$ ! 
Every hypothesis is verified for Waldhausen's approximation theorem, 
thus $K({\call B}_{n+1}^{red}, {\call B}_n^{red}) \simeq K(\L)
\simeq K(\C_R,homology \; equivalences) \simeq K(R)$. \\
$\bullet$ Let's now some extra-precautions when dealing with 
``reduced objects'' : first, remark from \cite{bil2} that if a 
$\cD_0$-compex $C_*$ is reduced, then testing if $C_*$ is 
${\call B}_n$-local or ${\call B}_n^{red}$-local is the same; 
moreover, we know also \cite{bil2} that $K({\call B}_{n+1},
{\call B}_n) \simeq K({\call B}_{n+1}^{red},{\call B}_n^{red})$.  
Now by our main theorem, this is the K-theory of reduced objects 
$C_*$ for which $C_i$ is contractible $\forall 0 \le i \le n$, 
and the other $C_i$ equal $C_{n+1} \in \C_R$ via the homology 
equivalences $\lambda_i$. One could object that the contractible 
object $E_*$ that we built is not reduced : that's true, but 
in \cite{bil2} we have a functorial way to send the global diagram 
to reduced objects : one compose on the right to obtain the wanted 
$C_* \to D_* \to D_*^{red}$. We finish by remarking that 
$D_{n+1}^{red} \leftarrow D_{n+1} \to D$ is better than two 
homology equivalences : as ``acyclic=contractible'' in $\C_R$, 
\mbox{here are two opposite homotopy equivalences ( they have 
some inverse ! ), so $D_{n+1}^{red} \simeq D$. $\blacksquare$} \\
\end{demo}

\no
\ding{183} Now, let's exercise our excision theorem on this 
same class ${\call B}$. \\
Let ${\call F \subset B}_{n+1}$ be the class of $\cD_0$-complexes 
$(0 \to B_1 = B_1 = \ldots )$. \\
Let ${\call G \subset B}_{n+1}$ be the class of $\cD_0$-complexes 
$(0 \to 0 \to B_2 \to B_3 \ldots )$. \\
Then we have the following homotopy equivalences of Waldhausen
categories : ${\call F} \simeq \C_R$ via the projection functor and 
${\call G} \simeq {\call B}_n$ via the translation functor.\\
Moreover ${\call F \cap G} = 0$ : there is no non-trivial map from 
${\call F}$ to ${\call G}$. \\ Hence the transversality : 
${\call F \pitchfork G}$, and we have ${\call F + G = B}_{n+1}$. \\
The excision theorem gives the equivalence of spectra : 
$K({\call B}_{n+1},{\call B}_n) \simeq K(R)$. By an obvious 
induction on $n$, we obtain the result : \hfill
\fbox{$K({\call B}_n) \simeq K(R)^n$}.  \\

\no
\ding{184} We shall now study the nilpotent case ${\call A}_n 
\subset {\call A}_{n+1}$.

\begin{prop}.\\
Let $C_*$ be a reduced $\cD_0$-complex in $\tilde {\call A}$. \\
Then ``$C_*$ is ${\call A}_n$-local'' if and only if ``the map 
$\lambda$ induces a homology equivalence between $Ker(\alpha_m)$ 
and $Ker(\alpha_{m+1})$ for all index $1 \le m \le n$''. We can 
visualize more easily the ${\call A}_n$-local multicomplexes by 
writing horizontally the maps $\lambda$ and vertically the maps 
$\alpha$, the symbol \textcircled{p} designing the 'exact' squares 
( ie squares that are both homotopy-pushout and homotopy-pullback ) :
$$\xymatrix@-1.88pc{
0 \ar@{ >->}[rr] \ar@{->>}[dd] && C_1 \ar@{ >->}[rr] \ar@{->>}[dd] 
&& C_2 \ar@{ >->}[rr] \ar@{->>}[dd] && C_3 \ar@{->>}[dd] & \ldots &
C_{n-1} \ar@{ >->}[rr] \ar@{->>}[dd] && C_n \ar@{ >->}[rr] 
\ar@{->>}[dd] && C_{n+1} \ar@{->>}[dd] & \ldots \\
&&& \mbox{\textcircled{p}} && \mbox{\textcircled{p}} &&&& 
\mbox{\textcircled{p}} \\
0 \ar@{ >->}[rr] && 0 \ar@{ >->}[rr] && C_1 \otimes S \ar@{ >->}[rr]
&& C_2 \otimes S & \ldots & C_{n-2} \otimes S \ar@{ >->}[rr] &&
C_{n-1} \otimes S \ar@{ >->}[rr] && C_n \otimes S & \ldots }$$ 
\end{prop}
\begin{demo}
$\bullet \underline{(i)\to(ii)}$ : Consider the test-objects $D_*=
g_m^{m+1}(D)=(0 \to \ldots 0 \to D \to CD = \ldots )$, 
beginning with $m$ zeros, where the complex $D=( \ldots \to 0 \to 
R \to 0 \to \ldots)$ is the base ring concentrated in degree $0$, 
and \mbox{$CD$ is the mapping-cone of the identity on $D$ ( that's
contractible ),} and every index $1 \le m < n$. Let $\cH om(D_*,
C_*)$ denote the set of 'algebraic morphisms' $f:D_* \to C_*$. By 
$R$-linearity, $f$ is given by the image $x=f(\bI_R) \in C_m$. 
To commute with $\alpha$, we must have $x \in Ker(\alpha_m)$. 
At the index $m+1$, the cone $CD$ is given by two copies of $R$ : 
one generated by $\bI_R$ in degree $0$, that is the image of $D$
by $\lambda$, and the other generated by $e$ in degree $-1$. 
The map $f_{m+1}$ is thus given by the images $\lambda(x)=
f_{m+1}(\bI_R) \in Ker(\alpha_m)$ and $y=f_{m+1}(e) \in 
Ker(\alpha_{m+1})$. It remains to complete on the right side by
the composed maps : $f_{m+1+k}=\lambda_{m+k} \circ \ldots 
\lambda_{m+1} \circ f_{m+1}$ to have the commutations with every
structural map. Hence we have an isomorphism of graded $R$-modules 
$\cH om(D_*,C_*) \simeq Ker(\alpha_m) \times Ker(\alpha_{m+1})$ 
defined by $f \leftrightarrow (x,y)$. Now describe the 
differential \vspace*{.05cm}through this correspondence : Leibniz' 
formula $[f \mapsto d \circ f -(-1)^{deg(f)}f\circ d]$\vspace*{.05cm}
applied to $D$ gives $[x \mapsto dx]$; on the cone $CD$, it gives 
$[y \mapsto dy -(-1)^{deg(f)}\lambda(x)]$; \vspace*{-.1cm}hence the 
correspondence $[df \leftrightarrow (dx,dy-(-1)^{deg(f)}\lambda(x)]$.
Next we consider the short exact sequence : 
$\xymatrix@-1pc{
Ker(\alpha_{m+1} ) _{p+1} \ar@{ >->}[r]^i & \cH om(D_*,C_* ) _p 
\ar@{->>}[r]^\pi & Ker(\alpha_m ) _p}$  
where the injection is given by $i(y)=(0,y)$ respecting 
$i(dy)=(0,dy)$ and the projection is given by $\pi(x,y)=x$ respecting
$\pi(d(x,y))=dx$. The two maps $i$ and $\pi$ are thus chain morphisms
and induce a long exact sequence in homology; moreover, $C_*$ being 
${\call A}_n$-local and $D_* \in {\call A}_n$, thus the complex 
$\cH om(D_*,C_*)$ is acyclic; hence we obtain the following 
isomorphism in homology $\delta:H_p(Ker(\alpha_m)) \simeq H_p(Ker(
\alpha_{m+1}))$. We must identify this connecting morphism : let 
$x \in Ker(\alpha_m)$ be a cycle, such that $dx=0$. We choose a 
section $(x,0) \in \cH om(D_*,C_*)$, then its boundary 
$d(x,0)=-(-1)^{deg(f)} \lambda(x)$. The connecting morphism 
thus corresponds ( up to a question of sign ) to $\lambda$. 
Conclusion : $\lambda$ induces a homology equivalence between 
$Ker(\alpha_m)$ and $Ker(\alpha_{m+1})$ for all index $1 \le m <n$.\\
$\bullet \underline{(ii) \to (i)}$ : The converse is much harder, 
and shall be proven directly for all index $n$ by the construction 
of an explicit homotopy, that is also available for the case 
$n=\infty$ ( in other words, it means the case of 
${\call A}$-local objects ). Let's start by reviewing our notations. 
Let $B_*$ be a reduced $\cD_0$-complex in ${\call A}$ ( that means 
$B_\infty$ is contractible ) of the type :       
$$\xymatrix@-1.8pc{
\ldots & B_{-1} \ar@{ >->}[rr] \ar@{->>}[dd] && B_0 
\ar@{ >->}[rr] \ar@{->>}[dd] && B_1 \ar@{ >->}[rr] \ar@{->>}[dd] 
&& B_2 \ar@{->>}[dd] & \ldots \\
&&&& \mbox{\textcircled{p}} && \mbox{\textcircled{p}} \\
\ldots & B_{-2} \otimes S \ar@{ >->}[rr] && B_{-1} \otimes S 
\ar@{ >->}[rr] && B_0 \otimes S \ar@{ >->}[rr] && B_1 \otimes S 
& \ldots }$$
where all squares for positive index are ``exact''. Consider a 
test-object $A_*=( \ldots 0 \to A_0 \to A_1 \to A_2 \ldots )$ in 
${\call A}$, doted with a map of $\cD_0$-complexes $f_*:A_* \to B_*$
given by a family of maps $\hat{f}_i:A_i \to B_i$ for all $i \ge 0$. 
We want to prove that $f \sim 0$. Take a family of compatible 
algebraic splittings for the following short exact sequences : 
( on the right side for $n \ge 0$ ) 
$$\xymatrix@-.5pc{
A_n \ar@<-2pt>@{ >->}[r]_{\lambda_n} & A_{n+1} 
\ar@<-2pt>@{->>}[l]_(.45){u_n} \ar@<-2pt>@{->>}[r]_(.45){\pi_{n+1}} 
& C_{n+1} \ar@<-2pt>@{ >->}[l]_{v_{n+1}} 
\ar `u^l[ll] `^dl[ll]_{\phi_{n+1}} [ll] &
& K \ar@<-2pt>@{ >->}[r]_{j_n} & B_n 
\ar@<-2pt>@{->>}[l]_(.45){\theta_n} 
\ar@<-2pt>@{->>}[r]_(.35){\beta_n} & 
B_{n-1} \otimes S \ar@<-2pt>@{ >->}[l]_(.6){\sigma_n}
\ar `u^l[ll] `^dl[ll]_{\delta_{n-1}} [ll]}$$
with the two structural maps :\hfill \fbox{\parbox{2.8cm}{
$\alpha_n:A_n \to A_{n-1} \otimes S \\ \mu_n:B_n \to B_{n+1}$}} \\
and the compatibility relation : \hfill
\fbox{$\theta_{n+1} \mu_n = \theta_n$} \\
and the differentials : \hspace{1.5cm} \fbox{\parbox{2.5cm}{
$d\theta_n = \delta_{n-1} \beta_n \\ 
d\sigma_n = -j_n \delta_{n-1} \\
du_n = \phi_{n+1} \pi_{n+1} \\
dv_{n+1} = - \lambda_n \phi_{n+1}$}} \hfill
\parbox{4cm}{( every other map is a chain morphism, in degree $0$ 
for $\lambda,\pi,j,\beta,\alpha,\mu$ and in degree $-1$ for 
the ``boundaries'' $\delta$ and $\phi$ )} \\

\no
We shall moreover need some specific notations ( independent 
from $n$ ) : \\

\no
$\forall p \ge 0$, we shall note : \hfill 
\fbox{\parbox{4.5cm}{$T_0 = \delta_n j_n \hspace{.92cm} : K \otimes S
\hspace{.42cm} \to K \\
T_p = \delta \underbrace{\sigma \sigma \ldots \sigma}_{p \; times}
j_n : K \otimes S^{p+1} \to  K$}} 

\begin{lemme}.\\
These chain morphisms in degree $-1$ verify the functional relation :
\hfill \fbox{$dT_p=\build{\sum}_{i+j=p-1}^{} T_i T_j$}. 
\end{lemme}
\begin{demo}
We already know that $d\delta_k = dj_k = 0$ for all $k$, thus the 
differential $dT_p$ has only terms with $d\sigma_k=-j_k\delta_{k-1}$.
By enumerating them ( $k^{th}$ between $p$ choices, from the end ),
and forgetting the $(. \otimes S)$ useless and cumbersome, 
we write :\\
$\begin{array}{lll}
dT_p &=& \build{\sum}_{k=1}^{p} (-1)\delta_{n+p} 
\sigma_{n+p} \ldots \sigma_{n+k+1}(-j_{n+k} \delta_{n+k-1}) 
\sigma_{n+k-1} \ldots \sigma_{n+1} j_n \\
&=& \build{\sum}_{k=1}^{p} (\delta_{n+p} \sigma_{n+p} \ldots 
\sigma_{n+k+1} j_{n+k} ) ( \delta_{n+k-1} \sigma_{n+k-1} \ldots
\sigma_{n+1} j_n )\\
&=&\build{\sum}_{k=1}^{p} T_{p-k} T_{k-1} = \build{\sum}_{i+j=p-1}^{}
T_i T_j \hfill \blacksquare \end{array}$ 
\end{demo}

\begin{lemme}.\\
We procede to different re-writings and re-interpretations of the 
maps coming from $f_*$ : we begin from $\hat{f}_n : A_n \to B_n$ 
first, that we decompose on the maps $f_n:A_n \to K$, and then in a 
more synthetic manner, we write them from a unique $F:A \to K$. The 
important intermediary formulas are the following : 
$$\fbox{$\hat{f}_{n+1}=j_{n+1} \theta_n \hat{f}_n u_n + \sigma_{n+1} 
\hat{f}_n u_{n+1} + j_{n+1} f_{n+1} \pi_{n+1}$}$$
$$\fbox{$f_n=F\lambda_n^\infty v_n$}$$
$$\fbox{$\hat{f}_n =\build{\sum}_{k=0}^{n} \sigma_n \ldots 
\sigma_{k+1} j_k F \lambda_k^\infty \alpha_{k+1} \ldots \alpha_n$}$$
\end{lemme}
\begin{demo}
$\bullet$ Compatibility of the maps of $\cD_0$-complexes : from the 
following commutative diagram, we shall write the maps $\hat{f}_n:A_n
\to B_n$ inductively on $n$, from the maps $f_n:A_n \to K$ ( we shall
procede by necessary, then sufficient condition ).\\
${} \hskip.5cm
\xymatrix{A_n \ar@{ >->}@<-2pt>_{\lambda_n}[r] \ar[d]_{\hat f_n} &
A_{n+1} \ar@{->>}@<-2pt>_{u_n}[l] \ar@{->>}_{\alpha_{n+1}}[r] 
\ar@{..>}[d]_{\hat f_{n+1}} & A_n \otimes S \ar^{\hat f_n \otimes 
\bI_S}[d]\\
B_n \ar@{ >->}_{\mu_n}[r] & B_{n+1} \ar@{->>}@<-2pt>_{\beta_{n+1}}[r]
& B_n \otimes S \ar@{ >->}@<-2pt>_{\sigma_{n+1}}[l]}$ \hskip1.3cm
\raisebox{-1.6pc}{\parbox{4cm}{We must then verify  
the two conditions by induction : 
\fbox{\parbox{3.8cm}{
$\hskip.7cm \hat f_{n+1} \lambda_n = \mu_n \hat f_n$ \\
${}\; \; \; \; \beta_{n+1} \hat f_{n+1}=\hat f_n \alpha_{n+1}$}}}}\\
The naive approach uses the algebraic sections, then corrects the 
eventual defect : let's pose $\hat f_{n+1} = \mu_n \hat f_n u_n + 
\sigma_{n+1} \hat f_n \alpha_{n+1} + \W$. We write the two 
conditions, using the induction hypothesis, and the two equalities :
$Id_{A_{n+1}}=\lambda_n u_n + v_{n+1} \pi_{n+1}$ and $Id_{B_{n+1}}=
j_{n+1} \theta_{n+1} + \sigma_{n+1} \beta_{n+1}$. We obtain :
$\W= -\sigma_{n+1} \mu_{n-1} \beta_n \hat f_n u_n + ( j_{n+1} 
\theta_{n+1} \W v_{n+1} \pi_{n+1} )$. The second mixed term 
can be interpreted as $j_{n+1} f_{n+1} \pi_{n+1}$ with some $f_{n+1}:
C_{n+1} \to K$ on which we have no condition. Grouping the terms 
finishing by $\hat f_n u_n$, we can write in a synthetic manner :
$$\hat f_{n+1} = j_{n+1} \theta_n \hat f_n u_n + \sigma_{n+1} \hat 
f_n u_{n+1} +j_{n+1} f_{n+1} \pi_{n+1}$$
$\bullet$ Now we shall think algebraically at the $f_n:C_n \to K$ 
as coordinates of a unique map $F:A \to K$ where the total space is 
$A=A_\infty=\oplus C_n$, and the maps $\lambda_n^\infty:A_n \to A$ 
give the projections on each factor : $f_n=F\lambda_n^\infty v_n$.
We have a correspondence : $(f_n) \leftrightarrow (F:A \to K)$.
So for a homotopy we shall search a map $(g_n:C_n \to K) 
\leftrightarrow ( \hat g_n= d \hat f_n )$ rather in the form 
$(g_n) \leftrightarrow (G:A \to K)$. \\
$\bullet$ Let's precise the conventions for the third formula : 
if $k=n$, then the compositions with $\sigma$ on the left and 
$\alpha$ on the right are both empty; if $k=n-1$, we have only 
one term $\sigma_n$ and $\alpha_n$; and if $k<n-1$, then the index
are going down on the left side, and up on the right side. 
Now, the method for finding this formula is 'magic' : try and guess
for small index, then prove it officially by induction. 
For $n=0$, we get : $\hat f_0=j_0 f_0 \pi_0=j_0 F \lambda_0^\infty 
v_0 \pi_0=j_0 F \lambda_0^\infty$ because $C_0=A_0$. Suppose now the 
formula available for the index $n$ : we want to calculate 
$\hat f_{n+1}$. We replace $\hat f_{n+1}$ by the first formula; 
$f_{n+1}$ by $F \lambda_{n+1}^\infty v_{n+1}$; and $v_{n+1} 
\pi_{n+1}$ by $Id- \lambda_n u_n$. We then obtain the wanted result, 
plus the defect : $j_{n+1} ( \theta_n \build{\sum}_{k=0}^{n} 
\sigma_n \ldots \sigma_{k+1} j_k F \lambda_k^\infty \alpha_{k+1} 
\ldots \alpha_n - F \lambda_n^\infty) u_n$. But the composition 
$\theta_n \sigma_n$ is trivial, hence every term vanishes except 
for $k=n$. It remains : $j_{n+1} ( \theta_n j_n F \lambda_n^\infty 
-F \lambda_n^\infty ) u_n$. But here $\theta_n j_n=Id$ and the 
parenthesis vanishes. Thus the third formula is proven by induction.
\hfill $\blacksquare$ 
\end{demo}
\begin{lemme}.\\
Differentiating the induction formulae for $\hat f_n$, we obtain 
the following formula, giving $G$ associated with $\hat g_n=d 
\hat f_n$ : \hfill
\fbox{$dF-G=\build{\sum}_{i \ge 0}^{} T_i F \alpha^{i+1}$}.\\
Then we build a new differential $\delta$ defined by :\hfill
\fbox{$F \mapsto \delta F=dF - \build{\sum}_{i \ge 0}^{} T_i F 
\alpha^{i+1}$}.\\
According to the principle for inverting formal series, we obtain 
for every cycle $\delta F=0$, the general formula giving a homotopy 
$\delta G=F$ :
$$\fbox{$G=\build{\sum}_{p \ge 0}^{} ( \build{\sum}_{i_1, \ldots,i_p}
^{} (-1)^{(p+1).deg(F)} T_{i_p} T_{i_{p-1}} \ldots T_{i_1} F k 
\alpha^{i_1 +1} k \alpha^{i_2 +1} k \ldots k \alpha^{i_p +1} k )$}$$
\end{lemme}
\begin{demo}
$\bullet$ We differentiate the formula giving $\hat f_{n+1}$, and 
identify $\hat g_n=d \hat f_n$. After simplification, we obtain :
$g_{n+1}=df_{n+1}+(-1)^{deg(f)} \theta_n \hat f_n \phi_{n+1} -
\delta_n \hat f_n \alpha_{n+1} v_{n+1}$. Then we differentiate the 
expression : $f_n=F \lambda_n^\infty v_n$ and we find the value of 
$G$ by identification, from its coordinates $g_{n+1}=G \lambda_{n+1}
^\infty v_{n+1}$ projected on $C_{n+1}$ : $g_{n+1}=dF \lambda_{n+1}
^\infty v_{n+1} - \build{\sum}_{k=0}^{n} T_{n-k} F \lambda_k ^\infty 
\alpha^{n-k} \alpha v_{n+1}$. In a more compact way, we proved by 
projection the following formula : 
$dF -G = \build{\sum}_{i \ge 0}^{} T_i F \alpha^{i+1}$. \\
$\bullet$ After re-interpretation, we have $\cH om(A_*,B_*) =
(\cH om(A,K),\delta : F \mapsto G)$ where the new differential 
is given by \mbox{$f \mapsto \delta F = dF - \build{\sum}
_{i \ge 0}^{} T_i F \alpha ^{i+1}$. Let's verify that $\delta$ is 
really a differential :} 
$$\delta^2 F=d^2 F - \build{\sum}_{i \ge 0}^{} (dT_i)F \alpha^{i+1} 
+ \build{\sum}_{i \ge 0}^{} T_i (dF) \alpha^{i+1} - \build{\sum}_{i
  \ge 0}^{} T_i ( dF - \build{\sum}_{j \ge 0}^{} T_j F \alpha^{j+1} )
\alpha^{i+1}$$        
The double differential $d^2 F$ vanishes, and the term $i=0$ from 
the first sum also, because $dT_0=0$; the second sum simplifies 
with the beginning of the third sum, and it remains : \\
$\delta^2 F= - \build{\sum}_{i \ ge 1}^{} ( \build{\sum}_{k+l=i-1}
^{} T_k T_l F \alpha^{i+1} ) + \build{\sum}_{i \ge 0}^{} ( \build
{\sum}_{j \ge 0}^{} T_i T_j F \alpha^{i+j+2})$; but after 
re-indexation, these two sums both equal $\build{\sum}
_{a \ge 0}^{} ( \build{\sum}_{k+l=a}^{} T_k T_l F \alpha^{a+2} )$,
thus $\delta^2 F=0$ and $\delta$ is really a differential ! \\ 
$\bullet$ We know that $A$ is acyclic : hence there exists a homotopy
$k:A \to A$ of degree $1$ such that $Id_A=dk$. Let $F$ be a cycle 
such that $dF=0$. Then $d(Fk)=(-1)^{deg(f)} F$ : it's a solution to 
``$F = \delta G$'' at the order $0$ in $T$ ! The problem is now, 
knowing that $d(k \alpha k)=\alpha k - k \alpha$, to guess a solution
to find our homotopy, according to the ( rather complicated ) 
principle for inverting formal series. Here I prefer to parachute 
the result ( found by P. Vogel ) coldly, and we shall verify it 
formally.  
$$G=\build{\sum}_{p \ge 0}^{} ( \build{\sum}_{i_1, \ldots,i_p}
^{} (-1)^{(p+1).deg(F)} T_{i_p} T_{i_{p-1}} \ldots T_{i_1} F k 
\alpha^{i_1 +1} k \alpha^{i_2 +1} k \ldots k \alpha^{i_p +1} k )$$
Let's write methodically its differential $\delta G$ and we shall 
obtain four triple sums : the first corresponds to juxtapose $T_i$ 
on the left, and $\alpha ^{i+1}$ on the right; the second to derive 
one of the $T_j$; the third is the differential $dF$ ( which we write
as a sum, from the hypothesis that $F$ is a cycle : $\delta F=0$ );
at last, the fourth corresponds to derive one of the $k$. \\
$\delta G= - \build{\sum}_i^{} ( \build{\sum}_p^{} ( \build{\sum}
_{i_1, \ldots, i_p}^{} (-1)^{(p+1).deg(F)} T_i T_{i_p} \ldots T_{i_1}
F k \alpha^{i_1 +1} k \ldots k \alpha^{i_p +1} k 
\alpha^{i+1})) \\
\hspace*{.8cm} +\build{\sum}_p^{} ( \build{\sum}_{i_1, \ldots,i_p}^{}
( \build{\sum}_{w=1}^p ( \build{\sum}_{l+m=i_w -1}^{} (-1)^{(p+1).
deg(F)} (-1)^{w-1} T_{i_p} \ldots (T_l T_m) \ldots T_{i_1} \\
\hspace*{8.3cm} F k \alpha^{i_1 +1} k \ldots k \alpha^{i_p +1} k)))\\
\hspace*{.8cm} +\build{\sum}_p^{} ( \build{\sum}_{i_1, \ldots,i_p}^{}
( \build{\sum}_i^{} (-1)^{(p+1).deg(F)} (-1)^p T_{i_p}\ldots T_{i_1} 
(T_i F \alpha^{i+1} ) k \alpha^{i_1 +1} k \ldots k \alpha^{i_p +1} 
k )) \\
\hspace*{.8cm} +\build{\sum}_p^{} ( \build{\sum}_{i_1, \ldots,i_p}^{}
( \build{\sum}_{w=0}^p (-1)^{(p+1).deg(F)} (-1)^{deg(F)+p+w} T_{i_p} 
\ldots T_{i_1} \\ \hspace*{5.8cm}
F k \alpha^{i_1 +1} k \ldots \alpha^{i_w +1} 
\alpha^{i_{w+1}+1} k \ldots k \alpha^{i_p +1} k ))$ \\

\no
Now re-interpret each of these sums : for the first, pose 
$i=i_{p+1}$, then it lacks a $k$ on the right; for the second, 
we cut $\alpha^{i_w +1}=\alpha^{l+1} \alpha^{m+1}$, then it lacks 
a $k$ between the two $\alpha$; for the third, pose $i=i_0$, then 
it lacks a $k$ on the left. With these simplifications, the fourth 
( with the index $p+1$ ) vanishes with the terms of the three others
( with the index $p$ ). For this, we take an index $1 \le w \le p$, 
and we impose three conditions for the annulation of the signs : 
( let $w'=p+1-w$ going down from $p$ to $1$ ), the conditions 
correspond respectively to the ``middle'', the ``beggining'' and 
the ``end'' for the lacking $k$.      
$$\fbox{\parbox{5.7cm}{
$\epsilon_{p+1}.(-1)^{deg(F)+p+1+w} + \epsilon_p . (-1)^{w' -1}=0 \\
\epsilon_{p+1} . (-1)^{deg(F)+p+1} + \epsilon_p . (-1)^p =0 \\
\epsilon_{p+1} . (-1)^{deg(F)+p+1+p+1}  - \epsilon_p =0$}}$$
Remark : the notation $\epsilon$ allows us to make the theoretical 
verification of the signs in the sums indexed by $(i_1,\ldots,i_p)$.
Actually, all three conditions are equivalent to : $\epsilon_{p+1} 
= \epsilon_p . (-1)^{deg(F)}$. Hence the existence of a sign : 
$\epsilon_p=\epsilon_0 . (-1)^{p.deg(F)}$. But we have seen in our 
first development 'by the hand' that $\epsilon_0=(-1)^{deg(F)}$. 
Thus the formula vanishes well. It then remains two minor points of 
detail : first, the change of index in the ``middle'' term. We 
enumerate the integral points in the plane by their usual coordinates
$(m,l)$ or by the diagonals $(l+m=i_w -1, w \ge 1)$ : actually, the 
case $i_w=0$ has a trivial differential. Finally, the compensation 
lets untouched the term $p=0$ in the fourth sum : here, we have 
$\epsilon_0 . (-1)^{deg(F)} . F =F$, so we have proved that 
$\delta G=F$. $\blacksquare$ 
\end{demo}
\hfill
This ends the proof of Proposition 5. $\blacksquare$ \\
\end{demo}

\no
{\bf Remark } : \\
We lack a good interpretation of these ${\call A}_n$-local 
objects, to apply Waldhausen's approximation theorem, and allow 
us an ersatz of Quillen's d\'evissage for the difficult class  
${\call A}$ of nilpotent $\cD_0$-complexes. 
Actually, if one manages to prove $K({\call A}_{n+1},{\call A}_n) 
\simeq K(R)$, then Vogel's excision [Conjecture] would be solved.  
\newpage
\nocite{*}
\bibliography{local}

\vskip 1cm
\no
{\bf I would like to thank my Master, 
Mr Pierre VOGEL, who gave me much of his time 
and many of his ideas and techniques, so 
that this article could be published. 
Thanks for his gentleness and his always wise 
counsels. I think he should really be listed 
as a co-author.}

\begin{center}
\vspace{1cm}
\begin{tabular}[t]{l}
{\bf Mr BIHLER Frank} \\
\\
U.F.R. de Math{\'e}matiques -- Case 7012 \\
Universit{\'e} de Paris 7 -- Denis Diderot \\
2, Place Jussieu \\
75251 PARIS cedex 05 \\
FRANCE \\
\\
{\bf  Mail :} bihler@math.jussieu.fr \\
{\bf Webpage :} http://www.institut.math.
jussieu.fr/{\~{}}bihler
\end{tabular} 
\end{center}

\end{document}